\pdfoutput=1
\documentclass[letterpaper, 10pt, journal, twoside]{IEEEtran}

\usepackage{amsmath,amssymb,amsthm,mathtools}
\usepackage{graphicx}
\usepackage{booktabs}
\usepackage{tikz}
\usepackage{pgfplots}
\pgfplotsset{compat=1.18}
\pgfplotsset{table/search path={data}}%
\usetikzlibrary{arrows.meta,positioning,calc}
\usepackage[hidelinks]{hyperref}
\usepackage{cite}

\definecolor{cA}{HTML}{0072B2}
\definecolor{cB}{HTML}{D55E00}
\definecolor{cC}{HTML}{009E73}
\definecolor{cD}{HTML}{666666}
\definecolor{cE}{HTML}{CC79A7}

\pgfplotsset{
  lcsty/.style={
    width=\linewidth, height=5.6cm, tick align=outside,
    tick label style={font=\scriptsize},
    label style={font=\scriptsize},
    title style={font=\scriptsize\bfseries, yshift=-2pt},
    grid=major, grid style={gray!20, line width=0.3pt},
    axis line style={gray!55, line width=0.5pt},
    every axis plot/.append style={line width=0.9pt},
    legend style={font=\tiny, draw=gray!45, fill=white, fill opacity=0.92,
                  text opacity=1, inner sep=1.5pt, row sep=-1.5pt},
    legend cell align=left,
  },
}


\newtheorem{theorem}{Theorem}[section]

\newtheorem{proposition}[theorem]{Proposition}
\newtheorem{lemma}[theorem]{Lemma}

\theoremstyle{definition}
\newtheorem{definition}[theorem]{Definition}
\newtheorem{assumption}[theorem]{Assumption}
\theoremstyle{remark}

\newtheorem{observation}[theorem]{Numerical Observation}

\newcommand{\E}{\mathbb{E}}
\newcommand{\R}{\mathbb{R}}
\newcommand{\Normal}{\mathcal{N}}

\newcommand{\given}{\,|\,}

\DeclareMathOperator{\Var}{Var}
\DeclareMathOperator{\Cov}{Cov}

\begin{document}

\title{\LARGE A Multistage Extension of Witsenhausen's\\ Counterexample with Linear Signaling}

\author{Quanyan Zhu and Tamer Ba\c{s}ar%
\thanks{Q. Zhu is with the Department of Electrical and Computer
Engineering, New York University, Brooklyn, NY 11201, USA
(e-mail: \texttt{quanyan.zhu@gmail.com}).}%
\thanks{T. Ba\c{s}ar is with the Coordinated Science Laboratory,
University of Illinois Urbana-Champaign, Urbana, IL 61801, USA.}}

\markboth{}%
{Zhu and Ba\c{s}ar: A Multistage Extension of Witsenhausen's Counterexample}

\maketitle

\begin{abstract}
We study a multistage scalar extension of Witsenhausen's counterexample in
which controller~1 observes the initial state exactly, whereas every
subsequent controller observes only a local noise-corrupted state and does
not inherit its predecessors' observations, so that every nonterminal action
both regulates the state and shapes the signal available to the next
controller. Restricting the signaling controllers to linear policies while
leaving the terminal controller unrestricted, we compute the exact value
within this class and reduce the optimization to a deterministic dynamic
program in a single propagated variance. The optimal intermediate gains are shown to have a shrinkage form, and the
value functions of that program are concave. For two stages, the benchmark admits a closed form on a matching locus, where
for small effort weights it has two distinct minimizers, and at the standard
parameter instance on that locus a two-point nonlinear policy strictly
improves upon it. For three stages, numerical continuation at a reference instance reveals
three distinct local minima, and two values of the initial variance at which
the lowest of the three changes abruptly. Under
the terminal-cost normalization we further prove, for every finite horizon,
that the optimal linear value becomes asymptotically equal to the
signaling-idle benchmark as the initial variance tends to infinity.
\end{abstract}


\section{Introduction}

\IEEEPARstart{W}{itsenhausen's} counterexample \cite{Wit68} showed that a
two-stage scalar linear-quadratic-Gaussian team can have a nonlinear optimal
policy once the information structure is not classical. The underlying
mechanism is by now well understood: the first controller may expend effort
not merely to regulate the state, but also to \emph{encode} it, so that the
second controller, which observes only a noise-corrupted version of what the
first has produced, is able to estimate it more accurately.

The corresponding problem over longer horizons is less completely
characterized. With three or more stages, an intermediate controller must decode
an upstream signal and, through one and the same action, re-encode it for its
successor. This multihop mechanism has no counterpart in the two-stage
problem, and it is the object of study of the present Letter.

Ho and Chu \cite{HoChu72} established affine optimality for partially nested
linear-quadratic-Gaussian teams, and Bansal and Ba\c{s}ar \cite{BanBas87} characterized conditions under
which affine optimality is retained; the problem studied here sits outside
that class. We establish existence for the present problem using Witsenhausen's static
reduction \cite{Wit88} together with the team results of
\cite{GYBL15,Yuk20}; see \cite{YukBas13} for information structures in
networked control.

The two-stage counterexample has been attacked numerically in several ways,
including approximating networks \cite{BPZ01}, hierarchical search
\cite{LLH01}, and methods carrying approximation guarantees \cite{GPS13}. Closest to our formulation, Wu and
Verd\'u \cite{WuVer11} recast it as a quadratic Wasserstein transport
problem coupled to a mean-square error term; that view does not carry over
unchanged, since only the first hop is free transport. A finite-level restriction has been developed in \cite{Zhu}. Closest in theme
among multistage studies, Lipsa and Martins \cite{LipMar11} have shown for a
scalar memoryless-control problem that a two-valued policy can beat the
affine optimum once the horizon is long enough; their controllers act on a
common noisy measurement of a Markov state rather than on the state each
predecessor produced, and thus the signal a controller receives is not shaped by
its predecessors' actions, and the multihop mechanism studied here does not
arise.

We obtain the exact value of the problem over the class in which the signaling
controllers are linear and the terminal controller is unrestricted, and show
that it collapses to a deterministic dynamic program in a single propagated
variance whose intermediate stationary gains are shrinkage factors in $[0,1)$. At $N=2$, the benchmark has a closed-form expression on the locus
$k^2\sigma_0^2=R$ and is attained at two gains; at $N=3$, a numerical study
exhibits three distinct local minima; and for every finite horizon we prove
that the
optimal linear value becomes asymptotically equal to the signaling-idle
benchmark as the initial variance grows.

Section~\ref{sec:formulation} sets out the model, its classical counterpart
and the resulting bracket on the team value; Section~\ref{sec:linear}
develops the benchmark and its two- and three-stage specializations;
Section~\ref{sec:conclusions} concludes.

\section{A Multistage Extension and Its Classical Counterpart}
\label{sec:formulation}

The system consists of $N$ controllers acting sequentially on a single scalar state, with quadratic penalties on the nonterminal control actions
and on the resulting states. What distinguishes it from a standard sequential problem is the
information pattern: the first controller observes the state exactly, every later one observes only a noise-corrupted version of what its predecessor had produced, and none inherits its predecessors' observations. We also record the
classical counterpart obtained by restoring nestedness, whose value brackets
the team cost from below.

\begin{assumption}\label{ass:prim}
Fix $N\ge2$. On a common probability space, the initial state $X_0$ and the
observation noises $V_1,\ldots,V_{N-1}$ are mutually independent, zero-mean
Gaussian random variables, $X_0\sim\Normal(0,P_0)$ and
$V_i\sim\Normal(0,R_i)$, with $P_0>0$ and $R_i>0$ for every $i$. The effort
weights satisfy $k_i>0$ for $i=1,\ldots,N-1$.
\end{assumption}

\begin{definition}[Multistage extension of the Witsenhausen
counterexample]\label{def:model}
Controller~$1$ observes $X_0$ and sets $U_1=\gamma_1(X_0)$, $X_1=X_0+U_1$; for $i=2,\ldots,N$, controller~$i$ observes $Y_i=X_{i-1}+V_{i-1}$ and sets
$U_i=\gamma_i(Y_i)$, $X_i=X_{i-1}+U_i$. The information fields are
\begin{equation}
  \mathcal I_1=\sigma(X_0),\ \ \mathcal I_i=\sigma(Y_i),\ \ i\ge2,
  \label{eq:fields}
\end{equation}
and $\Gamma$ denotes the set of Borel policy tuples
$\gamma=(\gamma_1,\ldots,\gamma_N)$ with $\E[U_i^2]<\infty$. With state
weights $c_i\ge0$ and $c_N>0$, the cost of $\gamma\in\Gamma$ is
\begin{equation}
  J(\gamma)=\E\Bigl[\textstyle\sum_{i=1}^{N-1}k_i^2U_i^2+\sum_{i=1}^{N}c_iX_i^2\Bigr],
  \ \ J^*=\inf_{\gamma\in\Gamma}J(\gamma).
  \label{eq:obj}
\end{equation}
\end{definition}

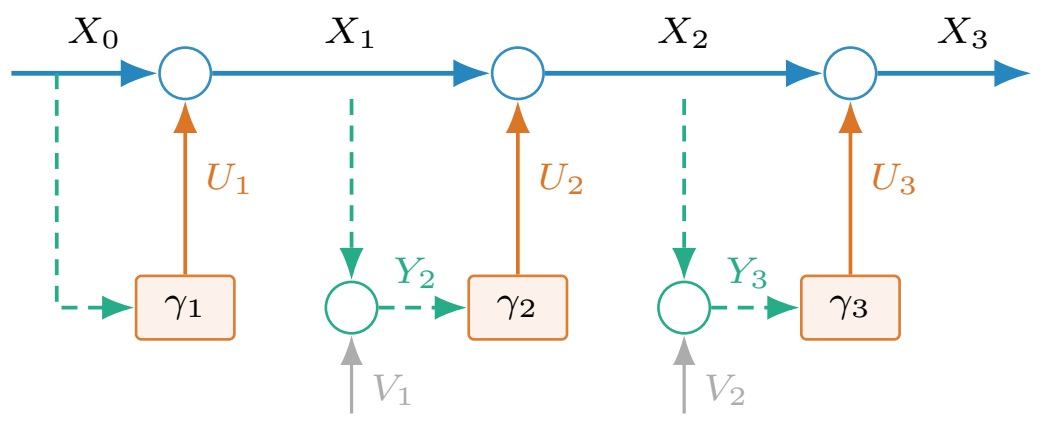
\begin{figure}[t]
\centering
\begin{tikzpicture}[>=Latex, font=\scriptsize,
  sum/.style={circle, draw=cA!80, line width=0.5pt, inner sep=0pt,
              minimum size=3.4mm},
  osum/.style={circle, draw=cC!85, line width=0.5pt, inner sep=0pt,
               minimum size=3.4mm},
  ctl/.style={draw=cB!80, fill=cB!8, line width=0.5pt, rounded corners=1pt,
              minimum width=6.5mm, minimum height=4.2mm},
  st/.style={cA!85, line width=0.9pt},
  inf/.style={cC!85, densely dashed, line width=0.7pt},
  act/.style={cB!85, line width=0.7pt},
  nz/.style={gray!65, line width=0.6pt}]

  \coordinate (x0) at (0,0);
  \node[sum] (s1) at (1.15,0) {};
  \node[sum] (s2) at (3.35,0) {};
  \node[sum] (s3) at (5.55,0) {};
  \coordinate (end) at (6.75,0);

  \draw[st,->] (x0) -- (s1);   \draw[st,->] (s1) -- (s2);
  \draw[st,->] (s2) -- (s3);   \draw[st,->] (s3) -- (end);
  \node[above=0.4mm] at (0.55,0) {$X_0$};
  \node[above=0.4mm] at (2.25,0) {$X_1$};
  \node[above=0.4mm] at (4.45,0) {$X_2$};
  \node[above=0.4mm] at (6.30,0) {$X_3$};

  \node[ctl] (c1) at (1.15,-1.55) {$\gamma_1$};
  \node[ctl] (c2) at (3.35,-1.55) {$\gamma_2$};
  \node[ctl] (c3) at (5.55,-1.55) {$\gamma_3$};
  \draw[act,->] (c1) -- node[right,pos=0.55] {$U_1$} (s1);
  \draw[act,->] (c2) -- node[right,pos=0.55] {$U_2$} (s2);
  \draw[act,->] (c3) -- node[right,pos=0.55] {$U_3$} (s3);

  \node[osum] (o2) at (2.25,-1.55) {};
  \node[osum] (o3) at (4.45,-1.55) {};
  \draw[inf,->] (0.30,0) -- (0.30,-1.55) -- (c1);
  \draw[inf,->] (2.25,-0.17) -- (o2);
  \draw[inf,->] (4.45,-0.17) -- (o3);
  \draw[inf,->] (o2) -- node[above,pos=0.42] {$Y_2$} (c2);
  \draw[inf,->] (o3) -- node[above,pos=0.42] {$Y_3$} (c3);
  \draw[nz,->] (2.25,-2.25) -- node[right,pos=0.30] {$V_1$} (o2);
  \draw[nz,->] (4.45,-2.25) -- node[right,pos=0.30] {$V_2$} (o3);
\end{tikzpicture}
\caption{The multistage system at $N=3$. Blue: the state path. Orange: the control
actions. Green dashed: what each controller sees. Gray: the noises.
Controller~$1$ alone sees the state exactly; later ones see only a
noise-corrupted version of what the predecessor produced, and do \emph{not}
inherit $Y_2,\ldots,Y_{i-1}$, so every nonterminal action is at once
regulator and encoder.}
\label{fig:arch}
\end{figure}

See Fig.~\ref{fig:arch}. Dividing \eqref{eq:obj} by $c_N$ leaves the
minimizers unchanged, so one may normalize $c_N=1$; we retain $c_N$ in the
general formulas below and impose the normalization in the two- and
three-stage specializations. The further choice $c_1=\cdots=c_{N-1}=0$, under which the intermediate
states are not charged for at all and the cost reduces to signaling effort
plus terminal error, we call the \emph{terminal-cost normalization}. It is
the setting of Witsenhausen's original two-stage problem, and by
Proposition~\ref{prop:classical} the one in which the classical counterpart
costs nothing. We omit an effort penalty on $U_N$, since the terminal action has no
downstream signaling role. Such a weight would replace the terminal
criterion by a positive combination of the terminal-state error and the
second moment of $U_N$; the linear-class analysis below can be adapted after
the terminal value function is modified accordingly. For $N=2$ and under the terminal-cost normalization, \eqref{eq:obj} is the
counterexample of \cite{Wit68}, up to the sign convention adopted for the
actions.

The structure \eqref{eq:fields} is dynamic, since $Y_i$ depends on
$U_1,\ldots,U_{i-1}$, and nonclassical: partial nestedness fails already at
$i=2$, where $U_1$ affects $Y_2$ yet $\sigma(Y_2)\not\supseteq\sigma(X_0)$,
and analogous local failures occur at every subsequent link. Replacing
\eqref{eq:fields} by its nested closure $\mathcal
I_i^{\mathrm{cl}}=\sigma(X_0,Y_2,\ldots,Y_i)$, with dynamics and cost untouched, yields a problem whose value
$J^{\mathrm{cl}}$ has a closed-form expression.

\begin{proposition}\label{prop:classical}
Set $\rho_N:=0$ and, for $i=N-1,\ldots,1$,
\begin{equation}
  \rho_i=\frac{k_i^2(c_i+\rho_{i+1})}{k_i^2+c_i+\rho_{i+1}} .
  \label{eq:rho}
\end{equation}
Then, $J^{\mathrm{cl}}=\rho_1P_0$, attained by the linear state feedback
$U_i=-(c_i+\rho_{i+1})X_{i-1}/(k_i^2+c_i+\rho_{i+1})$ for $i<N$ and
$U_N=-X_{N-1}$, and $J^{\mathrm{cl}}=0$ under the terminal-cost normalization. Moreover
\begin{equation}
  J^{\mathrm{cl}}\le J^*\le J_L^*\le J_{\mathrm{ol}},
  \label{eq:bracket}
\end{equation}
with $J_L^*$ the linear benchmark of Definition~\ref{def:linear} and
$J_{\mathrm{ol}}$ the \emph{signaling-idle} benchmark, obtained by setting
$a_1=\cdots=a_{N-1}=0$ and optimizing controller~$N$.
\end{proposition}

\begin{IEEEproof}
By induction $X_{i-1}$ is $\mathcal I_i^{\mathrm{cl}}$-measurable: this
holds at $i=1$, and if $X_{j-1}$ is, then so are $U_j$ and
$X_j=X_{j-1}+U_j$, hence also $\mathcal I_{j+1}^{\mathrm{cl}}$-measurable by
nestedness. Hence controller~$i$ can reconstruct the current state exactly, and since
the dynamics carry no process noise the whole trajectory is a deterministic
function of $X_0$ and the controls; the observation noises therefore carry
no additional payoff-relevant information about the state. Conditionally on $X_0=x$, the problem is thus the deterministic program
\[
  \min_{u_1,\ldots,u_N}\ \sum_{i=1}^{N-1}k_i^2u_i^2+\sum_{i=1}^{N}c_ix_i^2,
  \qquad x_i=x_{i-1}+u_i,\ \ x_0=x .
\]
This is a scalar linear-quadratic regulator, solved backwards in the
standard way \cite{BasOls99}: the terminal stage contributes
$\min_uc_N(x_{N-1}+u)^2=0$ at $u_N=-x_{N-1}$, so $\rho_N=0$; and if the
cost-to-go from $x_i$ is $\rho_{i+1}x_i^2$, then stage $i$ minimizes
$k_i^2u^2+q(x_{i-1}+u)^2$ with $q:=c_i+\rho_{i+1}\ge0$, whose value is
$k_i^2q\,x_{i-1}^2/(k_i^2+q)$ at $u=-qx_{i-1}/(k_i^2+q)$ since $k_i^2+q>0$.
This is \eqref{eq:rho} with the stated feedback. Expectation over $X_0$ gives
$J^{\mathrm{cl}}=\rho_1P_0$; and if $c_i=0$ for $i<N$ then $\rho_i=0$ by
downward induction.

For \eqref{eq:bracket}, $\mathcal I_i\subseteq\mathcal I_i^{\mathrm{cl}}$
for every $i$, and thus every tuple admissible for the problem is admissible for its
counterpart and the infimum over the larger set is no larger; the other two
inequalities hold because $\Gamma_L\subseteq\Gamma$ and the signaling-idle
tuple lies in $\Gamma_L$.
\end{IEEEproof}

All other bounds below are \emph{upper} bounds on $J^*$, and thus \eqref{eq:bracket} furnishes the only lower one; under the terminal-cost normalization it is
vacuous, which is the sense in which that normalization is the hardest case for the
lower bound in \eqref{eq:bracket}.

Before restricting the policy class, we record that the unrestricted problem is
well posed.

\begin{theorem}[Existence]\label{thm:exist}
Under Assumption~\ref{ass:prim} there exists a deterministic Borel policy
tuple $\gamma^\star\in\Gamma$ with $J(\gamma^\star)=J^*$.
\end{theorem}

\begin{IEEEproof}
For $i\ge2$, the conditional law of $Y_i$ given $X_0$ and $U_1,\ldots,U_{i-1}$ has the Lebesgue density $\varphi_{R_{i-1}}\bigl(y-X_0-\sum_{j<i}U_j\bigr)$, where, for a variance $R>0$, $\varphi_R(y):=(2\pi R)^{-1/2}\exp(-y^2/2R)$ denotes the density of $\Normal(0,R)$. Since $\varphi_{R_{i-1}}(y)>0$ for every $y\in\R$, this law is mutually absolutely continuous with respect to the fixed reference law
$\Normal(0,R_{i-1})$, which does not depend on the actions. It is absolute
continuity, and not pointwise domination, that the reduction requires: the
ratio $\varphi_R(y-\mu)/\varphi_R(y)=\exp(\mu y/R-\mu^2/2R)$ is unbounded in
$y$ for $\mu\neq0$. Witsenhausen's static reduction \cite{Wit88} therefore
applies, and transferring the action dependence from the measure to the
integrand yields an equivalent static team whose measurements are mutually
independent and policy independent, with the fields \eqref{eq:fields}, the
policy space $\Gamma$ and the value $J^*$ unchanged.

The problem is thus a dynamic team admitting a static reduction with independent
measurements. For that class, continuity of the reduced cost in the actions
suffices for a team-optimal policy to exist, and the optimizer may be taken
deterministic, provided the induced strategic measures are tight
\cite{Yuk20}. Continuity holds because the reduced cost is a polynomial in the
actions multiplied by Gaussian likelihood ratios. For tightness, restrict to
the sublevel set $\{\gamma\in\Gamma:J(\gamma)\le J_{\mathrm{ol}}\}$, which is
nonempty and contains a minimizing sequence; on it $k_i^2\E[U_i^2]\le
J_{\mathrm{ol}}$ for each $i<N$. This is a bound under the original measure,
and converting it into tightness of the strategic measures under the reference
measure is exactly the step carried out in \cite[\S5]{GYBL15} for $N=2$, by
means of the function class of \cite[Def.~4.4]{GYBL15} and
\cite[Lem.~4.5]{GYBL15}; each stage cost has the same form, so
the argument applies at each stage in turn. Transporting the optimizer back
through the reduction gives $\gamma^\star$.
\end{IEEEproof}

\section{The Linear Signaling Class}
\label{sec:linear}

Since no closed-form solution of the unrestricted problem is known even at
$N=2$, we fix a tractable class and compute its value exactly: the signaling
controllers, whose actions carry information forward, are restricted to be
linear, while the terminal controller, which merely estimates, is left
unrestricted. The resulting upper bound on $J^*$ collapses to a scalar
dynamic program whose minimizing gains are confined to an interval.

\begin{definition}\label{def:linear}
The linear signaling class $\Gamma_L\subset\Gamma$ consists of the policy
tuples of the form
\begin{equation}
  U_1=-a_1X_0,\ \ U_i=-a_iY_i,\ \ 2\le i\le N-1,
  \label{eq:linpol}
\end{equation}
with $a_1,\ldots,a_{N-1}\in\R$ and $\gamma_N$ unrestricted. We write
$J_L^*=\inf_{\gamma\in\Gamma_L}J(\gamma)$.
\end{definition}

It is only the \emph{signaling} controls that are restricted here; the
terminal controller is left unrestricted, and when the signaling policies are linear, every state is Gaussian and the optimal terminal controller is therefore linear. This is well known for Witsenhausen's two-stage problem \cite{Wit68}; the following lemma records its straightforward multistage version, included for completeness.

\begin{lemma}\label{lem:rec}
Let $\gamma\in\Gamma_L$ and set $P_i:=\Var(X_i)$. Then, every $X_i$ is zero-mean Gaussian,
\begin{equation}
  P_1=(1-a_1)^2P_0,\ \ P_i=(1-a_i)^2P_{i-1}+a_i^2R_{i-1},
  \label{eq:var}
\end{equation}
and $\E[U_1^2]=a_1^2P_0$, $\E[U_i^2]=a_i^2(P_{i-1}+R_{i-1})$ for
$2\le i\le N-1$. Furthermore, among all Borel $\gamma_N$ with
$\E[\gamma_N(Y_N)^2]<\infty$ the terminal cost $c_N\E[(X_{N-1}+\gamma_N(Y_N))^2]$
is minimized uniquely, up to a null set, by
\begin{equation}
  U_N^*=-\E[X_{N-1}\given Y_N]=-\kappa_NY_N,\ \
  \kappa_N=\tfrac{P_{N-1}}{P_{N-1}+R_{N-1}},
  \label{eq:term}
\end{equation}
the optimal value being $c_NP_{N-1}R_{N-1}/(P_{N-1}+R_{N-1})$.
\end{lemma}

\begin{IEEEproof}
From $U_1=-a_1X_0$, we obtain $X_1=(1-a_1)X_0$, which is zero-mean Gaussian of
variance $(1-a_1)^2P_0$. For $2\le i\le N-1$, substituting $Y_i=X_{i-1}+V_{i-1}$
into $U_i=-a_iY_i$ and thence into $X_i=X_{i-1}+U_i$ gives
\begin{equation}
  X_i=(1-a_i)X_{i-1}-a_iV_{i-1}.
  \label{eq:xrec}
\end{equation}
By induction $X_{i-1}$ is a measurable function of $X_0,V_1,\ldots,V_{i-2}$
alone, hence zero-mean Gaussian and independent of $V_{i-1}$; both
properties pass to $X_i$ through \eqref{eq:xrec}, and independence gives
\eqref{eq:var}. Also $\E[U_1^2]=a_1^2P_0$, while for $i\ge2$ that same
independence gives $\Var(Y_i)=P_{i-1}+R_{i-1}$ with $Y_i$ zero-mean, so
$\E[U_i^2]=a_i^2(P_{i-1}+R_{i-1})$.

For the terminal stage, the conditional mean minimizes the mean-square error
among all square-integrable estimators, by the orthogonality principle
\cite{Ber17}; note that this uses only square integrability, not
Gaussianity, a point we shall need in Proposition~\ref{prop:twopoint}. Hence
$\gamma_N(Y_N)=-\E[X_{N-1}\given Y_N]$ almost surely, uniquely up to a null
set. By the first part $(X_{N-1},Y_N)$ is jointly Gaussian and zero-mean with
$\Var(Y_N)=P_{N-1}+R_{N-1}>0$ and $\Cov(X_{N-1},Y_N)=P_{N-1}$, so
$\E[X_{N-1}\given Y_N]=\kappa_NY_N$ and the conditional variance is
$P_{N-1}-\kappa_NP_{N-1}=P_{N-1}R_{N-1}/(P_{N-1}+R_{N-1})$. The factor
$c_N>0$ does not affect the minimizer.
\end{IEEEproof}

As far as downstream estimation is concerned, \eqref{eq:var} compresses the
entire upstream policy history into the single scalar $P_i$.

\begin{proposition}[Exact value over $\Gamma_L$]\label{prop:benchmark}
Under Assumption~\ref{ass:prim}, $J_L^*=\min_{a\in\R^{N-1}}J_L(a)$, where
\begin{equation}
\begin{aligned}
  J_L=\ &k_1^2a_1^2P_0+\sum_{i=2}^{N-1}k_i^2a_i^2(P_{i-1}+R_{i-1})\\
        &+\sum_{i=1}^{N-1}c_iP_i+\frac{c_NP_{N-1}R_{N-1}}{P_{N-1}+R_{N-1}},
\end{aligned}
  \label{eq:JL}
\end{equation}
with $P_1,\ldots,P_{N-1}$ generated by \eqref{eq:var}. The minimum is
attained.
\end{proposition}

\begin{IEEEproof}
Fix $a\in\R^{N-1}$. By Lemma~\ref{lem:rec} the control terms of \eqref{eq:obj}
are $k_1^2a_1^2P_0$ and $k_i^2a_i^2(P_{i-1}+R_{i-1})$, and the state terms are
$c_i\E[X_i^2]=c_iP_i$ for $i<N$; none of these involves $\gamma_N$. The one
remaining term is $c_N\E[X_N^2]=c_N\E[(X_{N-1}+\gamma_N(Y_N))^2]$, which
Lemma~\ref{lem:rec} minimizes over $\gamma_N$, giving \eqref{eq:JL}. Taking
the infimum over $a$ then gives $J_L^*$.

For attainment, $J_L$ is smooth, since each $P_i$ is a polynomial in
$(a_1,\ldots,a_i)$ and the terminal term is rational with denominator $\ge
R_{N-1}>0$. Setting $\alpha_1:=k_1^2P_0$ and $\alpha_i:=k_i^2R_{i-1}$ for
$2\le i\le N-1$, all positive, and $\kappa:=\min_{1\le i\le N-1}\alpha_i>0$,
nonnegativity of every summand together with $P_{i-1}+R_{i-1}\ge R_{i-1}$
gives $J_L\ge\sum_{i=1}^{N-1}\alpha_ia_i^2\ge\kappa\|a\|^2$. Thus, $J_L$ is coercive and attains its infimum. Indexing the coefficients this way avoids
an empty minimum when $N=2$.
\end{IEEEproof}

\subsection{A scalar dynamic program}

Expression \eqref{eq:JL} is a deterministic optimization in $N-1$ variables,
and its structure invites a backward recursion. Define
$F_N(P):=c_NPR_{N-1}/(P+R_{N-1})$ and, recursively for $i=N-1,\ldots,2$,
\begin{equation}
  F_i(P)=\min_{a\in\R}\bigl\{k_i^2a^2(P+R_{i-1})+c_iQ+F_{i+1}(Q)\bigr\},
  \label{eq:dp}
\end{equation}
where $Q=Q_i(a,P):=(1-a)^2P+a^2R_{i-1}$; the subscript on $Q_i$ records
that the noise variance appearing in it is the one entering stage~$i$.

\begin{lemma}\label{lem:F}
For $i=2,\ldots,N-1$, the minimum in \eqref{eq:dp} is attained. Moreover each $F_i:[0,\infty)\to[0,\infty)$, $i=2,\ldots,N$, is continuous,
nondecreasing, \emph{concave}, and satisfies $F_i(0)=0$.
\end{lemma}

\begin{IEEEproof}
Backward induction. For $i=N$ all assertions are immediate from the displayed
formula for $F_N$. Assume them for $F_{i+1}$ and set
$G(a,P):=k_i^2a^2(P+R_{i-1})+c_iQ_i(a,P)+F_{i+1}(Q_i(a,P))$, which is continuous
and nonnegative on $\R\times[0,\infty)$. Fix $\bar P>0$. Since
$G(a,P)\ge k_i^2R_{i-1}a^2$ uniformly in $P\in[0,\bar P]$, and $G(0,\cdot)$ is
bounded on $[0,\bar P]$, there is $A=A(\bar P)$ with $G(a,P)>G(0,P)$ whenever
$|a|>A$; hence the infimum over $a\in\R$ equals the minimum over the compact
set $\{|a|\le A\}$ and is attained. Continuity of $F_i$ on $[0,\bar P]$ then
follows from continuity of $G$ on the compact product
$[-A,A]\times[0,\bar P]$ by Berge's maximum theorem, and $\bar P$ was arbitrary.
For monotonicity, let $0\le P\le P'$; for each fixed $a$ both
$k_i^2a^2(P+R_{i-1})$ and $Q_i(a,P)$ are nondecreasing in $P$, and $c_i\ge0$, so
$G(a,P)\le G(a,P')$ by the inductive monotonicity of $F_{i+1}$, and taking
minima preserves the inequality. Finally, $0\le F_i(0)\le G(0,0)=c_i\cdot0+F_{i+1}(0)=0$. For concavity,
$F_N(P)=c_N\bigl(R_{N-1}-R_{N-1}^2/(P+R_{N-1})\bigr)$ is concave; and if
$F_{i+1}$ is concave, then, for each fixed $a$, the map $P\mapsto G(a,P)$ is
the sum of two affine functions of $P$ and of $F_{i+1}$ composed with the
affine map $P\mapsto Q_i(a,P)$, hence concave. A pointwise infimum of
concave functions is concave, and thus $F_i$ is concave.
\end{IEEEproof}

\begin{theorem}[Scalar dynamic program]\label{thm:dp}
Under Assumption~\ref{ass:prim},
\begin{equation}
  J_L^*=\min_{a_1\in\R}\bigl\{k_1^2a_1^2P_0+c_1P_1+F_2(P_1)\bigr\},
  \ \ P_1=(1-a_1)^2P_0 .
  \label{eq:dpouter}
\end{equation}
\end{theorem}

\begin{IEEEproof}
By Proposition~\ref{prop:benchmark} and the definition of $F_N$, minimizing
$J_L$ over $a\in\R^{N-1}$ subject to \eqref{eq:var} is a deterministic
finite-horizon optimal control problem: the state is $P_i\in[0,\infty)$, the
control is $a_i\in\R$, the stage cost $k_i^2a_i^2(P_{i-1}+R_{i-1})+c_iP_i$ and
the successor state $P_i=(1-a_i)^2P_{i-1}+a_i^2R_{i-1}$ depend on the past
only through the pair $(P_{i-1},a_i)$, and the terminal cost is
$F_N(P_{N-1})$. Recursion \eqref{eq:dp} is the associated Bellman equation,
and its minima are attained by Lemma~\ref{lem:F}; the dynamic programming
principle for deterministic finite-horizon problems
\cite[Ch.~1]{Ber17} therefore identifies $F_j(P)$, for $j=2,\ldots,N$, as the
least cost of running stages $j$ through $N$ from a variance $P$ entering
stage~$j$. Taking $j=2$ with $P=P_1=(1-a_1)^2P_0$ and minimizing over $a_1$
gives \eqref{eq:dpouter}, the outer minimum being attained by
Proposition~\ref{prop:benchmark}.
\end{IEEEproof}

Within $\Gamma_L$ the transmitted variance is thus a sufficient summary of
the upstream history; for nonlinear signaling it is the
\emph{shape} of the induced law, not merely its second moment, that the
successor exploits.

\begin{proposition}[Shrinkage]\label{prop:shrink}
Let $2\le i\le N-1$, let $P_{i-1}\ge0$, and let $a_i$ be \emph{any} minimizer
in \eqref{eq:dp} at $P=P_{i-1}$. Then
\begin{equation}
  0\ \le\ a_i\ \le\ \frac{P_{i-1}}{P_{i-1}+R_{i-1}}\ <\ 1 ,
  \label{eq:shrink}
\end{equation}
with no differentiability of $F_{i+1}$ assumed: no minimizing gain
overshoots, and none completely cancels the inherited state component. If, in addition, $F_{i+1}$ is differentiable at
$Q_i:=(1-a_i)^2P_{i-1}+a_i^2R_{i-1}$, then, writing
$\Lambda_i:=c_i+F_{i+1}'(Q_i)\ge0$ for the effective downstream price of
variance,
\begin{equation}
  a_i=\frac{\Lambda_iP_{i-1}}{(k_i^2+\Lambda_i)(P_{i-1}+R_{i-1})} .
  \label{eq:gain}
\end{equation}
\end{proposition}

\begin{IEEEproof}
If $P_{i-1}=0$, then both $Q_i(a,0)=a^2R_{i-1}$ and the effort term $k_i^2a^2R_{i-1}$ are strictly increasing in $|a|$, and thus $a_i=0$ is the unique minimizer. Assume
therefore $P_{i-1}>0$ and set $m:=P_{i-1}/(P_{i-1}+R_{i-1})\in(0,1)$. As a
function of $a$,
$Q_i(a,P_{i-1})=P_{i-1}-2aP_{i-1}+a^2(P_{i-1}+R_{i-1})$ is a quadratic with
minimum at $a=m$, hence symmetric about $m$. If $a<0$, then $Q_i(a,P_{i-1})>Q_i(0,P_{i-1})$, and thus the terms $c_iQ_i$ and $F_{i+1}(Q_i)$ are
no smaller at $a$ than at $0$, since $c_i\ge0$ and $F_{i+1}$ is nondecreasing
by Lemma~\ref{lem:F}, while the effort term is strictly positive at $a$ and
vanishes at $0$; thus no minimizer is negative. If $a>m$, put $a':=2m-a$;
symmetry gives $Q_i(a',P_{i-1})=Q_i(a,P_{i-1})$, leaving the state and
continuation terms unchanged, whereas $|a'|<|a|$ strictly lowers the effort
term, so no minimizer exceeds $m$. Finally, $m<1$ because $R_{i-1}>0$.

For \eqref{eq:gain}, $F_{i+1}'(Q_i)\ge0$ by Lemma~\ref{lem:F} and $c_i\ge0$, and thus $\Lambda_i\ge0$. Differentiating the bracket in \eqref{eq:dp} and using
$\partial_aQ_i=2[a(P_{i-1}+R_{i-1})-P_{i-1}]$, the terms $c_iQ_i$ and
$F_{i+1}(Q_i)$ together contribute $\Lambda_i\,\partial_aQ_i$, and thus half the derivative is $k_i^2a(P_{i-1}+R_{i-1})+\Lambda_i[a(P_{i-1}+R_{i-1})-P_{i-1}]$.
This vanishes at a minimizer, and $k_i^2+\Lambda_i\ge k_i^2>0$ permits
solving for $a_i$.
\end{IEEEproof}

This is a tradeoff between two prices: a larger downstream price $\Lambda_i$
induces controller~$i$ to transmit less variance, a larger local weight
$k_i^2$ suppresses its action. Since $c_i$ enters only as an increment to
$\Lambda_i$, charging for $X_i$ is indistinguishable, at the level of the
stationarity conditions, from raising the value the successor attaches to a
smaller variance. The weights $c_i$ therefore act as an additional penalty
on propagated state variance, that is, as a tax on carrying variance
forward.

\begin{theorem}[Asymptotic equivalence to signaling idle]\label{thm:futility}
Fix $N\ge2$, the noise variances $R_1,\ldots,R_{N-1}$ and the weights
$k_1,\ldots,k_{N-1}$, and assume the terminal-cost normalization. With
$J_{\mathrm{ol}}(P_0)=c_NP_0R_{N-1}/(P_0+R_{N-1})$, the signaling-idle cost,
\begin{equation}
  J_L^*(P_0)\longrightarrow c_NR_{N-1}
  \quad\text{and}\quad
  \frac{J_L^*(P_0)}{J_{\mathrm{ol}}(P_0)}\longrightarrow1
  \label{eq:futility}
\end{equation}
as $P_0\to\infty$. In words: once the initial variance is large, the best
linear policy costs essentially what it costs to set every gain $a_i$ to
zero and leave the whole job to the terminal estimator. Linear signaling
then buys nothing.
\end{theorem}

\begin{IEEEproof}
Normalize $c_N=1$. The upper bound is immediate, since
$J_L^*\le J_{\mathrm{ol}}<R_{N-1}$ and $J_{\mathrm{ol}}\to R_{N-1}$. For the
lower bound fix $\varepsilon\in(0,R_{N-1})$ and suppose
$J_L(a)\le R_{N-1}-\varepsilon$ for some $a\in\R^{N-1}$. Under the normalization, every term of \eqref{eq:JL} is nonnegative, hence bounded by
$R_{N-1}-\varepsilon$. The terminal term gives
\[
  P_{N-1}\ \le\ M:=\frac{R_{N-1}(R_{N-1}-\varepsilon)}{\varepsilon},
\]
and for $2\le i\le N-1$, the $i$th effort term gives
$a_i^2P_{i-1}\le a_i^2(P_{i-1}+R_{i-1})\le B_i:=(R_{N-1}-\varepsilon)/k_i^2$.
Since $P_i\ge(1-a_i)^2P_{i-1}$ by \eqref{eq:var} and
$2(1-a_i)^2+2a_i^2=2-4a_i+4a_i^2\ge1$,
\[
  P_{i-1}\ \le\ 2(1-a_i)^2P_{i-1}+2a_i^2P_{i-1}\ \le\ 2P_i+2B_i .
\]
Applying this for $i=N-1,N-2,\ldots,2$ in turn and starting from
$P_{N-1}\le M$, bounds $P_1$ by a constant $M_1$ depending only on
$\varepsilon$, on $N$, and on the $R_i$ and $k_i$, but not on $P_0$; for
$N=2$ the recursion is vacuous and $M_1=M$. As $P_1=(1-a_1)^2P_0$, this
forces $|1-a_1|\le\sqrt{M_1/P_0}$, hence $a_1\ge1-\sqrt{M_1/P_0}$, so that
once $P_0>M_1$,
\[
  R_{N-1}-\varepsilon\ \ge\ J_L(a)\ \ge\ k_1^2a_1^2P_0
  \ \ge\ k_1^2\bigl(\sqrt{P_0}-\sqrt{M_1}\bigr)^2 .
\]
The right-hand side diverges, so beyond some $\bar P(\varepsilon)$ no $a$
satisfies $J_L(a)\le R_{N-1}-\varepsilon$; that is,
$J_L^*(P_0)>R_{N-1}-\varepsilon$. As $\varepsilon$ was arbitrary,
$\liminf J_L^*\ge R_{N-1}$, and the ratio statement follows.
\end{IEEEproof}

\subsection{Two stages}

We now specialize to $N=2$, imposing the normalization $c_2=1$ and writing
$k:=k_1$, $\sigma_0^2:=P_0$ and $R:=R_1$, so that by
Proposition~\ref{prop:benchmark}
\begin{equation}
  J_L^{(2)}(a)=k^2a^2\sigma_0^2+c_1(1-a)^2\sigma_0^2
  +\frac{(1-a)^2\sigma_0^2R}{(1-a)^2\sigma_0^2+R}.
  \label{eq:JL2}
\end{equation}
On one distinguished parameter locus the minimum has a closed-form
expression and is attained at two distinct gains.

\begin{theorem}[Closed form on the matching locus]\label{thm:two}
Suppose $c_1=0$ and the matching condition $k^2\sigma_0^2=R$ holds. If
$0<k<\tfrac12$, then $J_L^{(2)}$ has exactly three stationary points; the
outer two,
\begin{equation}
  a^*_\pm=\tfrac12\bigl(1\mp\sqrt{1-4k^2}\bigr),\ \ a^*_-+a^*_+=1,
  \label{eq:twomin}
\end{equation}
are global minimizers with the common value
\begin{equation}
  J_L^{(2)*}=R\,(1-k^2),
  \label{eq:twoval}
\end{equation}
and the third, which lies strictly between them, is a local maximum. If
$k\ge\tfrac12$, the stationary point is unique and is the global minimizer.
\end{theorem}

\begin{IEEEproof}
Put $t:=1-a$, $u:=t^2\sigma_0^2$, and use $\frac{uR}{u+R}=R-\frac{R^2}{u+R}$.
Under $R=k^2\sigma_0^2$, one has $u+R=\sigma_0^2(t^2+k^2)$ and
$R^2/(u+R)=Rk^2/(t^2+k^2)$, so $J_L^{(2)}=R\varphi(t)$ with
\begin{equation}
  \varphi(t)=(1-t)^2+1-\frac{k^2}{t^2+k^2}.
  \label{eq:phi}
\end{equation}
Differentiating, $\varphi'(t)=-2(1-t)+2k^2t/(t^2+k^2)^2=-2q(t)/(t^2+k^2)^2$,
where
\[
  q(t):=(1-t)(t^2+k^2)^2-k^2t=-(t^2-t+k^2)(t^3+k^2t-k^2),
\]
the factorization being verified by expanding both sides. The stationary points of $\varphi$ are thus the roots of $q$; for
$0<k<\tfrac12$, the three real roots are simple.
The quadratic factor has roots $t_\pm=\tfrac12(1\pm\sqrt{1-4k^2})$, real and
distinct precisely when $0<k<\tfrac12$ and then satisfying
$0<t_-<\tfrac12<t_+<1$; the cubic $\psi(t):=t^3+k^2t-k^2$ is strictly
increasing with $\psi(0)<0<\psi(1)$, and thus it has one real root $t_0\in(0,1)$.

To locate them, use $k^2=t_\pm-t_\pm^2$ to compute
$\psi(t_\pm)=t_\pm^3+(t_\pm-t_\pm^2)t_\pm-(t_\pm-t_\pm^2)=t_\pm(2t_\pm-1)$,
which is negative at $t_-$ and positive at $t_+$; since $\psi$ is increasing,
$t_-<t_0<t_+$. Evaluating $\varphi$ at $t_\pm$, where $t^2+k^2=t$ and hence
$k^2/(t^2+k^2)=k^2/t=1-t$, gives
$\varphi(t_\pm)=(1-t_\pm)^2+1-(1-t_\pm)=1-t_\pm+t_\pm^2=1-k^2$, which is
\eqref{eq:twoval} after multiplication by $R$; and $a^*_\pm=1-t_\pm$ gives
\eqref{eq:twomin}, the two labellings being exchanged by $t_-+t_+=1$.

Finally $\varphi(t)\to\infty$ as $|t|\to\infty$, and thus $\varphi'<0$ near
$-\infty$ and $>0$ near $+\infty$; having three simple zeros, $\varphi'$
changes sign at each, giving the pattern $-,+,-,+$. The stationary points
are therefore, in increasing order, a minimum, a maximum and a minimum, and
the two outer ones share the value $R(1-k^2)$.

For the last assertion, if $k>\tfrac12$ the quadratic factor has no real root
and $t_0$ is the only stationary point, while if $k=\tfrac12$ then
$t_\pm=\tfrac12$ and $\psi(\tfrac12)=0$, so that all three coincide. A smooth
coercive function with a unique stationary point attains its global minimum
there.
\end{IEEEproof}

The matching condition equates the effort scale $k\sigma_0$ with the
observation-noise standard deviation $\sqrt R$ concealing the state. It is a
particularly unfavorable regime for linear policies, and the comparison with
doing nothing is sharp: on this locus $\sigma_0^2=R/k^2$, so the
signaling-idle cost is
$J_{\mathrm{ol}}=\sigma_0^2R/(\sigma_0^2+R)=R/(1+k^2)$, and dividing
\eqref{eq:twoval} by it gives
\begin{equation}
  \frac{J_L^{(2)*}}{J_{\mathrm{ol}}}=\frac{R(1-k^2)}{R/(1+k^2)}=1-k^4 .
  \label{eq:k4}
\end{equation}
Linear signaling therefore improves upon
leaving controller~1 idle by the fraction $k^4$ only: at $k=0.2$ by
$0.16\%$, which is why $0.96154$ falls merely to $0.96$. That the class is
not only this ineffective but strictly suboptimal follows from the standard
two-point construction below.

\begin{proposition}[Two-point policies]\label{prop:twopoint}
Let $N=2$ with $c_2=1$, write $J^{(2)*}$ for the unrestricted two-stage
value, and for $h>0$ let $\gamma_1^h(x):=h\operatorname{sgn}(x)-x$, so
that $X_1=h\operatorname{sgn}(X_0)$, and let $\gamma_2$ be the corresponding
optimal terminal policy. Then
\begin{equation}
  J_{2\mathrm{pt}}(h)=k^2\Bigl(h^2-2h\sigma_0\sqrt{\tfrac2\pi}+\sigma_0^2\Bigr)
  +c_1h^2+h^2\bigl(1-\Xi\bigr),
  \label{eq:twopoint}
\end{equation}
where $\Xi:=\E[\tanh^2(h(h+V_1)/R)]$, and consequently
$J^{(2)*}\le\inf_{h>0}J_{2\mathrm{pt}}(h)$.
\end{proposition}

\begin{IEEEproof}
For the control term,
$\E[U_1^2]=\E[(h\operatorname{sgn}(X_0)-X_0)^2]=h^2-2h\E|X_0|+\sigma_0^2$ with
$\E|X_0|=\sigma_0\sqrt{2/\pi}$, and the state term is
$c_1\E[X_1^2]=c_1h^2$ since $X_1^2\equiv h^2$. For the terminal term, write
$S:=\operatorname{sgn}(X_0)$, which is uniform on $\{\pm1\}$ and independent of
$V_1$, so that $X_1=hS$ and $Y_2=hS+V_1$. The likelihood ratio of $Y_2=y$
under $S=\pm1$ gives $\Pr(S=1\given Y_2=y)=e^{hy/R}/(e^{hy/R}+e^{-hy/R})$,
whence $\E[X_1\given Y_2=y]=h\tanh(hy/R)$. By the orthogonality principle, as noted in Lemma~\ref{lem:rec}, the optimal
terminal policy is $U_2=-\E[X_1\given Y_2]$ for this prior too, and
\[
  \E[X_2^2]=\E\bigl[\Var(X_1\given Y_2)\bigr]
  =\E\bigl[h^2-h^2\tanh^2(hY_2/R)\bigr],
\]
where $X_1^2\equiv h^2$ was used. Since $\tanh^2$ is even and $V_1$ is symmetric, averaging over $S=\pm1$
gives the same expectation as conditioning on $S=1$, and thus $Y_2$ may be
replaced by $h+V_1$ inside the expectation, which gives \eqref{eq:twopoint}. The final assertion holds because the pair
$(\gamma_1^h,\gamma_2)$ is admissible.
\end{IEEEproof}

At the standard instance $\sigma_0^2=25$, $R=1$, $k=0.2$, $c_1=0$, which
lies on the matching locus, Theorem~\ref{thm:two} gives $J_L^{(2)*}=0.96$. A
strict separation follows without any quadrature. Take $h=4$ in Proposition~\ref{prop:twopoint} and, in place of the optimal
terminal policy, use the hard decoder $U_2=-4\operatorname{sgn}(Y_2)$; the
resulting cost is an upper bound on $J_{2\mathrm{pt}}(4)$ and hence on
$J^{(2)*}$. Its first-stage cost is
$k^2\bigl(41-40\sqrt{2/\pi}\bigr)<0.3634$, and $X_2$ vanishes unless the
sign is decoded wrongly, so that
$\E[X_2^2]=64\,\Pr(V_1<-4)\le32e^{-8}<0.0108$ by the Gaussian tail bound
$\Pr(V_1<-x)\le\tfrac12e^{-x^2/2}$. Hence
\begin{equation}
        J^{(2)*}\ <\ 0.3742\ <\ 0.96\ =\ J_L^{(2)*},
  \label{eq:gap2}
\end{equation}
which recovers Witsenhausen's nonlinear-versus-linear separation
\cite{Wit68}. Numerically, optimizing $h$ in \eqref{eq:twopoint} by
quadrature lowers the left-hand side to about $0.3648$.

\subsection{Three stages}

We turn to $N=3$, the shortest horizon exhibiting the multihop
mechanism. Write $\sigma_1^2:=R_1$ and $\sigma_2^2:=R_2$, impose the normalization
$c_3=1$, and restrict to $U_1=-aX_0$ and $U_2=-bY_2$, so that
$P_1=(1-a)^2\sigma_0^2$ and $P_2=(1-b)^2P_1+b^2\sigma_1^2$, and hence
\begin{equation}
  J_L(a,b)=k_1^2a^2\sigma_0^2+c_1P_1+k_2^2b^2(P_1+\sigma_1^2)+c_2P_2
           +\tfrac{P_2\sigma_2^2}{P_2+\sigma_2^2}.
  \label{eq:JL3}
\end{equation}
The remainder of this subsection is numerical. Fix the reference instance
\begin{equation}
    \sigma_0^2=25,\ \ \sigma_1^2=\sigma_2^2=1,\ \ k_1=k_2=k=0.2,\ \ c_1=c_2=0,
  \label{eq:ref}
\end{equation}
so that the first-hop effort and noise scales satisfy the matching relation
$k^2\sigma_0^2=\sigma_1^2$ of Theorem~\ref{thm:two}; the second hop sees
source variance $P_1=(1-a)^2\sigma_0^2$ rather than $\sigma_0^2$, and so
meets that relation only for particular $a$. Two values frame the discussion: the signaling-idle cost is
$J_{\mathrm{ol}}=0.96154$, while cancelling the state outright at stage~1
costs $k^2\sigma_0^2=1$, so a useful linear policy must land below $0.96154$
and the margin is small. All reported local minima below were computed by multistart from a
uniform $11\times11$ grid over $[-0.2,1.2]^2$, polished by BFGS to gradient
tolerance $10^{-12}$, with values quoted to the digits stable under an
independent Nelder--Mead polish.

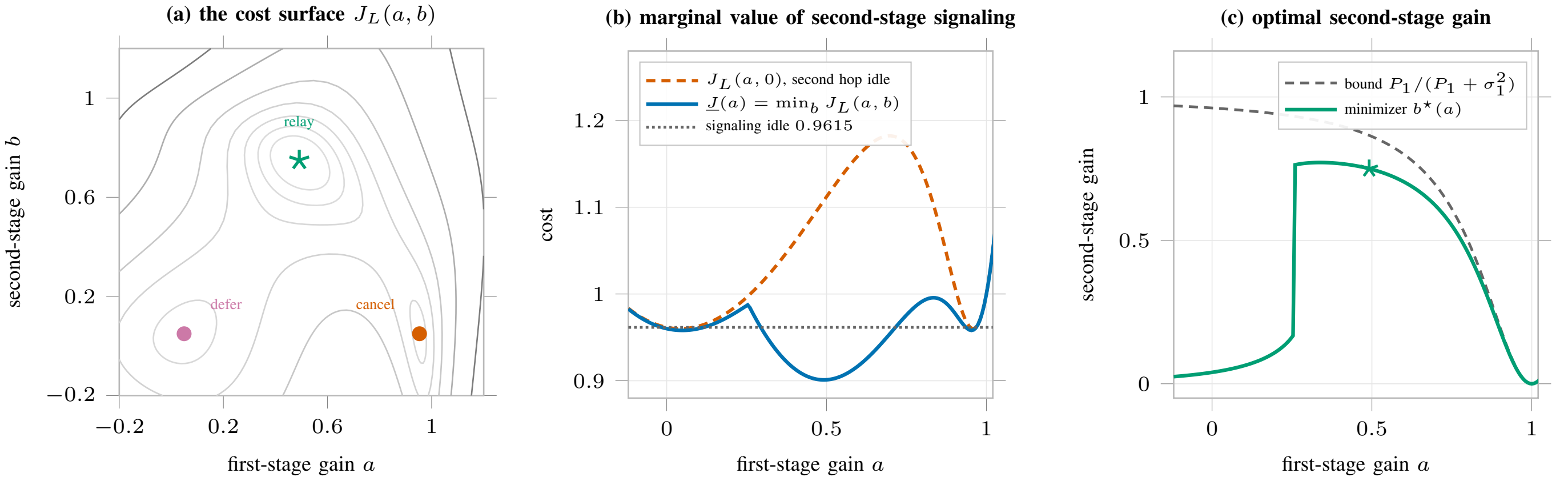
\begin{figure*}[!t]
\centering
\begin{minipage}[t]{0.32\textwidth}\centering
\begin{tikzpicture}
\begin{axis}[lcsty, xlabel={first-stage gain $a$},
  ylabel={second-stage gain $b$},
  title={(a) the cost surface $J_L(a,b)$},
  xmin=-0.2, xmax=1.2, ymin=-0.2, ymax=1.2,
  xtick={-0.2,0.2,...,1.2}, ytick={-0.2,0.2,...,1.2}, grid=none,
  colormap={lc}{color=(cD!22) color=(cD!85)}]
  \addplot[contour prepared, contour prepared format=standard,
           contour/labels=false, line width=0.5pt] table {fig-contour.dat};
  \addplot[only marks, mark=star, mark size=3.4pt, cC, line width=1pt]
    coordinates {(0.491485,0.750018)};
  \addplot[only marks, mark=*, mark size=2pt, cB]
    coordinates {(0.953329,0.049454)};
  \addplot[only marks, mark=*, mark size=2pt, cE]
    coordinates {(0.049110,0.049454)};
  \node[font=\tiny, text=cC, anchor=south] at (axis cs:0.4915,0.83) {relay};
  \node[font=\tiny, text=cB, anchor=east] at (axis cs:0.90,0.17) {cancel};
  \node[font=\tiny, text=cE, anchor=west] at (axis cs:0.11,0.17) {defer};
\end{axis}
\end{tikzpicture}
\end{minipage}\hfill
\begin{minipage}[t]{0.32\textwidth}\centering
\begin{tikzpicture}
\begin{axis}[lcsty, xlabel={first-stage gain $a$}, ylabel={cost},
  title={(b) marginal value of second-stage signaling},
  xmin=-0.12, xmax=1.02, ymin=0.88, ymax=1.28,
  legend pos=north west, legend columns=1]
  \addplot[cB, densely dashed, line width=1.1pt]
    table[x index=0, y index=3] {fig-reduced.dat};
    \addlegendentry{$J_L(a,0)$, second hop idle}
  \addplot[cA, line width=1.2pt]
    table[x index=0, y index=1] {fig-reduced.dat};
    \addlegendentry{$\underline J(a)=\min_bJ_L(a,b)$}
  \addplot[cD, densely dotted, line width=0.9pt]
    coordinates {(-0.12,0.961538) (1.02,0.961538)};
    \addlegendentry{signaling idle $0.9615$}
\end{axis}
\end{tikzpicture}
\end{minipage}\hfill
\begin{minipage}[t]{0.32\textwidth}\centering
\begin{tikzpicture}
\begin{axis}[lcsty, xlabel={first-stage gain $a$},
  ylabel={second-stage gain},
  title={(c) optimal second-stage gain},
  xmin=-0.12, xmax=1.02, ymin=-0.05, ymax=1.16,
  legend pos=north east, legend columns=1]
  \addplot[cD, densely dashed, line width=0.9pt] table {fig-bound.dat};
    \addlegendentry{bound $P_1/(P_1+\sigma_1^2)$}
  \addplot[cC, line width=1.2pt]
    table[x index=0, y index=2] {fig-reduced.dat};
    \addlegendentry{minimizer $b^\star(a)$}
  \addplot[only marks, mark=star, mark size=3pt, cC, forget plot]
    coordinates {(0.491485,0.750018)};
\end{axis}
\end{tikzpicture}
\end{minipage}
\caption{The three-stage linear benchmark at the reference instance
\eqref{eq:ref}. (a)~Level sets of $J_L(a,b)$ at
$J_L\in\{0.92,0.94,0.97,1.02,1.12,1.35,1.8\}$, darker being lower, with the
three local minima of Table~\ref{tab:three} marked. (b)~The partially
minimized cost $\underline J(a)$ against $J_L(a,0)$, the cost incurred with
controller~2 held idle; the vertical separation is the value of the second
signaling hop. (c)~The minimizing gain $b^\star(a)$ against the bound of
Proposition~\ref{prop:shrink}, which it respects everywhere; the absolute gap closes as $a\to1$ only
because both quantities vanish, the ratio tending to $0.9615$. The star marks the relay point.}
\label{fig:landscape}
\end{figure*}

\begin{figure*}[!t]
\centering
\begin{minipage}[t]{0.32\textwidth}\centering
\begin{tikzpicture}
\begin{axis}[lcsty, xlabel={$\sigma_0^2$}, ylabel={cost},
  title={(a) the three local minima and $J_L^*$},
  xmin=1, xmax=60, ymin=0, ymax=1.58,
  legend pos=south east, legend columns=2]
  \fill[cA!8] (axis cs:21.877,0) rectangle (axis cs:28.206,1.58);
  \addplot[cB, densely dashed] table {fig-br-cancel.dat};
    \addlegendentry{cancel}
  \addplot[cC, densely dotted, line width=1.1pt] table {fig-br-relay.dat};
    \addlegendentry{relay}
  \addplot[cE, dash dot] table {fig-br-defer.dat};
    \addlegendentry{defer}
  \addplot[cD, densely dashed] table[x index=0, y index=7] {fig-global.dat};
    \addlegendentry{signaling idle}
  \addplot[cA, line width=1.4pt] table[x index=0, y index=1] {fig-global.dat};
    \addlegendentry{$J_L^*=\min$}
  \addlegendimage{area legend, fill=cA!14, draw=cA!35}
    \addlegendentry{relay window}
\end{axis}
\end{tikzpicture}
\end{minipage}\hfill
\begin{minipage}[t]{0.32\textwidth}\centering
\begin{tikzpicture}
\begin{axis}[lcsty, xlabel={$\sigma_0^2$}, ylabel={optimal gain},
  title={(b) the gains jump at each threshold},
  xmin=1, xmax=60, ymin=-0.05, ymax=1.20,
  legend pos=north east, legend columns=2]
  \fill[cA!8] (axis cs:21.877,-0.05) rectangle (axis cs:28.206,1.20);
  \addplot[cA, line width=1.1pt] table[x index=0, y index=2] {fig-global.dat};
    \addlegendentry{$a^\star$}
  \addplot[cB, densely dashed, line width=1.1pt]
    table[x index=0, y index=3] {fig-global.dat};
    \addlegendentry{$b^\star$}
\end{axis}
\end{tikzpicture}
\end{minipage}\hfill
\begin{minipage}[t]{0.32\textwidth}\centering
\begin{tikzpicture}
\begin{axis}[lcsty, xlabel={$\sigma_0^2$}, ylabel={contribution to $J_L^*$},
  title={(c) cost decomposition},
  xmin=1, xmax=60, ymin=-0.04, ymax=1.08,
  legend style={at={(0.97,0.45)}, anchor=east}, legend columns=1]
  \fill[cA!8] (axis cs:21.877,-0.04) rectangle (axis cs:28.206,1.08);
  \addplot[cA, line width=1.1pt] table[x index=0, y index=4] {fig-global.dat};
    \addlegendentry{$k^2\E U_1^2$}
  \addplot[cC, dash dot, line width=1.1pt]
    table[x index=0, y index=6] {fig-global.dat};
    \addlegendentry{$\E[(X_3^\star)^2]$}
  \addplot[cB, densely dashed, line width=1.4pt]
    table[x index=0, y index=5] {fig-global.dat};
    \addlegendentry{$k^2\E U_2^2$}
\end{axis}
\end{tikzpicture}
\end{minipage}
\caption{Dependence on the initial variance at $\sigma_1^2=\sigma_2^2=1$,
$k=0.2$; the shaded band, with numerically determined boundaries at
approximately $21.88$ and $28.21$, is the window in which the two-hop relay
is the lowest-cost one found. (a)~The three local minima traced separately,
with $J_L^*$ their lower envelope; each persists well outside the range over
which it is the lowest, which motivates the use of multistart optimization when evaluating
$J_L^*$. (b)~The gains are discontinuous at both edges.
(c)~The three terms of \eqref{eq:JL3} that survive at $c_1=c_2=0$, plotted
separately rather than stacked: controller~2 expends appreciable effort only
inside the window, its contribution remaining below $1.5\times10^{-3}$
outside it and rising to $0.196$ within.}
\label{fig:sweep}
\end{figure*}
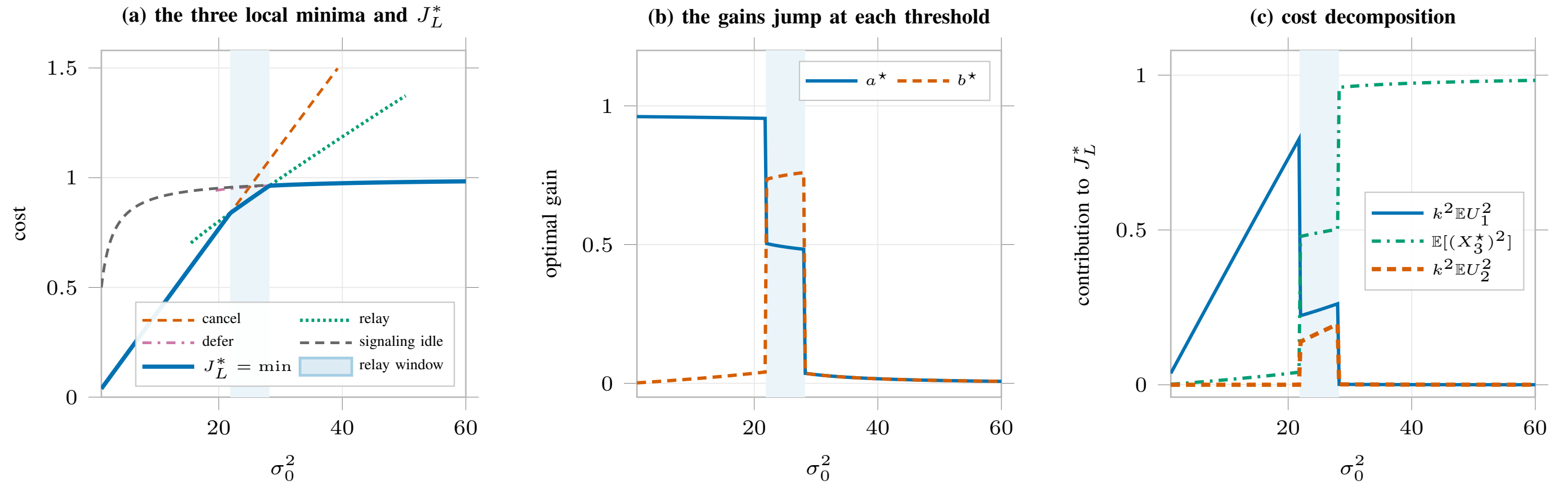

\begin{table}[!t]
\centering
\caption{The three local minima of $J_L$ at the reference instance
\eqref{eq:ref}. Relay, the only one of the three using both signaling actions appreciably,
is the
lowest-cost minimum found.}
\label{tab:three}
\scriptsize
\setlength{\tabcolsep}{3.2pt}
\begin{tabular}{@{}lrrrrrr@{}}
\toprule
 & $a^\star$ & $b^\star$ & $k^2\E U_1^2$ & $k^2\E U_2^2$
 & $\E[(X_3^\star)^2]$ & $J_L$\\
\midrule
Relay (lowest found) & 0.4915 & 0.7500 & 0.2416 & 0.1680 & 0.4915 & 0.9010\\
Cancel         & 0.9533 & 0.0495 & 0.9088 & 0.0001 & 0.0491 & 0.9581\\
Defer          & 0.0491 & 0.0495 & 0.0024 & 0.0023 & 0.9533 & 0.9581\\
\bottomrule
\end{tabular}
\end{table}

\begin{observation}[Three local minima]\label{obs:three}
At \eqref{eq:ref} our multistart computation consistently identifies three
local minima of $J_L$ (Table~\ref{tab:three}, Fig.~\ref{fig:landscape}(a)). \emph{Cancel} incurs almost all of its cost at stage~1, \emph{defer}
expends very little at either stage, and \emph{relay} divides the work
between the hops and is the lowest-cost minimum found, improving upon the
other two
by $5.95\%$ and upon the signaling-idle benchmark by $6.3\%$. All three respect Proposition~\ref{prop:shrink}, the relay bound reading
$b^\star\le0.8660$. The objective values of the two nonrelay minima agree to every computed
digit, a coincidence that recurs to eight digits at
$(k,\sigma_0^2)=(0.25,16)$ and $(0.1,100)$ and breaks in our computations
when $k\sigma_0\neq1$.
Fig.~\ref{fig:landscape}(b) sets $J_L(a,0)$, the cost with controller~2 held
idle, against $\underline J(a)=\min_bJ_L(a,b)$; the separation is the value
of the second hop. It attains $0.2439$ near $a=0.62$, to the right of the relay point, where it
is $0.2066$, and is numerically negligible
outside the central range, falling to $1.5\times10^{-3}$ at $a=0$ rather
than to zero. Panel~(c) shows the same in the minimizing gain: $b^\star(a)$ stays below $0.12$ for
$a\le0.2$, rises abruptly to $0.77$ by $a=0.4$, and returns to zero as
$a\to1$, respecting Proposition~\ref{prop:shrink} throughout. Its ratio to that
bound is $0.865$ at the relay point; by \eqref{eq:gain} the ratio
equals $\Lambda_2/(k^2+\Lambda_2)$, which tends to $1/(1+k^2)=0.9615$ as
$a\to1$, so the bound stays slack in relative terms even as both quantities
vanish.
\end{observation}

\begin{observation}[Two changes in the lowest minimum]\label{obs:regimes}
Sweeping $\sigma_0^2$ (Fig.~\ref{fig:sweep}(a),(b)), which of the three is
lowest changes twice, numerically at $\sigma_0^2\approx21.88$ and
$\sigma_0^2\approx28.21$: cancel is lowest below the first threshold, relay
between them, and defer above the second, that is $0.936<k\sigma_0<1.062$
for the relay. Both changes are abrupt: the two minima involved swap rank while remaining
separated in $(a,b)$, so the gains jump rather than drift. The relay optimum is also sensitive to the other parameters: at
$\sigma_1^2=2$ it survives only as a local minimum, its cost rising to
$0.9986$, and it ceases to be lowest at $c_1=0.010563$, two orders of magnitude below
$c_3=1$, so even a small intermediate-state penalty removes relay as the
lowest-cost linear solution in this computation.
\end{observation}

\section{Conclusions}
\label{sec:conclusions}
The linear signaling class in a scalar multistage system with nonclassical
information
admits an exact value, computable by a deterministic dynamic program in a
single propagated variance, with intermediate stationary gains confined to
$[0,1)$. At two stages, that value has a closed-form expression on the
matching locus and is attained at two distinct gains; at the standard
instance a two-point policy strictly improves upon it. At three stages, a
numerical study of a reference instance exhibits three distinct local
minima, with two abrupt changes in which of them is lowest, and we prove
that the optimal
linear value is asymptotically equal to the signaling-idle benchmark as the
initial variance grows. The benchmark thus reduces evaluation at any finite horizon to a
scalar-state dynamic program, and provides a tractable means of identifying
the parameter regimes in which multihop linear signaling is active.

Exhibiting any nonlinear policy of cost below $J_L^*$ would establish that
linear signaling is suboptimal at that horizon; quantifying how close such a
construction comes to the unrestricted optimum is a separate question, and
would require a lower bound sharper than \eqref{eq:bracket}. This multiplicity of local minima is consistent with the concavity of the
value functions
established in Lemma~\ref{lem:F}: minimizing a concave function of the
propagated variance is exactly the situation in which isolated local minima
and discontinuous exchanges between them arise.

\end{document}